\documentclass[a4paper,11pt]{article}
\usepackage{amssymb, amsmath, amsthm, amscd, amsfonts}

\newtheorem{theorem}{Theorem}

\newtheorem{example}{Example}
\newtheorem{algorithm}{Algorithm}
\newenvironment{pr}[1][Proof]{\noindent\textbf{#1.} }{\ \rule{0.5em}{0.5em}}

\begin{document}
\title{\textbf{Technical Notes on \lq\lq A new approach to the {\it d}-MC problem"}}
\author{\textbf{Majid~Forghani-elahabad\footnote{Email
address:~forghanimajid@mehr.sharif.ir} ,\
 Nezam~Mahdavi-Amiri\footnote{Corresponding author: Fax:++982166005117, Phone: ++982166165607, Email
address: nezamm@sina.sharif.edu} }}
\date{}
\maketitle $\vspace{-14mm}$
\begin{center}
\textit{{\scriptsize Faculty of Mathematical Sciences, Sharif University of 
Technology, Tehran, Iran}}\\
\end{center}

\hrule \section*{ Abstract} $\indent$ 
System reliability is the probability of the maximum flow in a stochastic-flow network from the source node to the sink node  being more than a demand level {\it d}. There are several approaches to compute system reliability using upper boundary points, called d-MinCuts ({\it d-MC}s). Search for all the {\it d-MC}s in a stochastic-flow network is an NP-hard problem. Here, a work proposed by Yeh [Yeh WC. A new approach to the d-MC problem. Reliab Eng and Syst Saf 2002; 77(2): 201-–206.] for determinig all the {\it d-MC}s is investigated. Two results (Lemma 3 and Theorem 5) are shown to be incorrect and their correct versions are established. Also, the complexity result (Theorem 6) is shown to be incorrect and the correct count is provided.

\vspace{2mm} \noindent {\normalsize{Keywords:}}
Reliability; Stochastic-flow network; d-MinCut ({\it d-MC}); Minimal Cuts ({\it MC}s).\\
\hrule
\section{ Introduction}

Reliability originates from a series of lectures given by Von Neumann in 1952 \cite{von}. After that, network reliability
theory has extensively been applied to a variety real-world systems such
as power transmission and distribution \cite{lin6},
computer and communication \cite{levitin}, transportation \cite{wu}, etc. Applying approximate methods \cite{coit2} or exact ones [6-10], system reliability can be computed in terms of the upper boundary points for demand level {\it d}, called d-MinCuts ({\it d-MC}s). Jan et al. \cite{janecc}, by introducing the notion of {\it d-MC} candidate, proposed an algorithm that first finds all the {\it d-MC} candidates obtained from each Minimal Cut ({\it MC}) and then checks every candidate for being a {\it d-MC}. Since then, most algorithms [7-10] have been composed of two general stages: finding all the {\it d-MC} candidates by an implicit enumeration method and verifying every candidate via a testing process for being a {\it d-MC}. Yeh \cite{yeh1} first presented an algorithm merely in accordance with the definition, which unfortunately had a defect \cite{salehi}. Proving some new results, Yeh \cite{yeh3} presented another algorithm. However, there are some flaws in \cite{yeh3} as will be discussed here. Yan and Qian \cite{yan} proved some new results to decrease the number of the obtained {\it d-MC} candidates, to find some {\it d-MC}s without the need for testing and to eliminate some duplicate {\it d-MC}s. Then, they proposed an improved algorithm which turned to be more efficient than the algorithms proposed in \cite{yeh1, yeh3} (see \cite{yan} for a comparative study). Salehi and Forghani \cite{salehi} proposed an algorithm by rectifying Yeh's algorithm in \cite{yeh1}. Here, by investigating the work proposed by Yeh \cite{yeh3}, we demonstrate that two results (Lemma 3 and Theorem 5) are incorrect and then, give and prove the correct versions of them. Moreover, the corresponding complexity result (Theorem 6) is shown to be incorrect and its correct version is presented.

In the remainder of our work, in Section 2 we provide the required definitions and illustrate the flaws of the proposed lemma and theorem using an example. Then, the correct versions are established. In Section 3, a computed time complexity of an existing algorithm is shown to be incorrect and its corresponding correct complexity result is established. Section 4 gives the concluding remarks.

\section{On the results}
Here, some required notations are described, and then two presented results (Lemma 3 and Theorem 5) in \cite{yeh3} are rewritten to explain their flaws and provide their correct forms.
\subsection{Problem description}
For convenience, the same notations, nomenclature, and assumptions used by Yeh \cite{yeh3} are given here. Let $G(V,E,W)$ be a stochastic-flow network with the set of nodes $V=\{1,2,...,n\}$, the set of edges $E=\{e_i| 1\leq i\leq m\}$, and $W(e_i)$ denoting the max-capacity of $e_i$, for $i=1,2,...,m$. The current capacity of arc $e_i$ is represent by $x_i$, and so $X = (x_1, x_2, ..., x_m)$ is a system-state vector representing the current capacity of all the arcs in $E$. Let $0(e_i)$ denote a system-state vector in which the capacity level is 1 for $e_i$ and 0 for other arcs, $C_i$ be the $i$th {\it MC} in $G(V,E,W)$, nodes $1$ and $n$ be the source and sink nodes, respectively, $p$ be the number of {\it MC}s in $G(V,E,W)$, $\sigma$ be the number of {\it d-MC} candidates obtained from each {\it MC}, $G(V,E,X)$ be the network corresponding to $G(V,E,W)$ with current state vector 
$X = (x_1, x_2, ..., x_m)$, $R(V,E,X^d)$ be the corresponding residual network to $G(V,E,X)$ after sending $d$ units of flow from node $1$ to node $n$, $W(X)$ be the max-flow from node $1$ to node $n$ in $G(V,E,X)$, $C^d_{ij}=(x_1, x_2, . . . , x_m)$ be the $j$th {\it d-MC} candidate (system-state
vector) generated from $C_i$ in $G(V, E, W)$, where $x_k = C^d_{ij}(e_k) \leq W(e_k),\ \sum x_k = d$, for $e_k \in C_i$ and $x_l = W(e_l),$ for $e_l \notin C_i$. Finally, a system-state vector $X$ is a {\it d-MC} if and only if $W(X)=d$ and $W(X+0(e_i))>d$, for every $e_i\in U(X)=\{e_i|X(e_i)<W(e_i)\}$.
\subsection{Incorrect results}

Two presented results (Lemma 3 and Theorem 5) in \cite{yeh3} have faults. Here,  we first rewrite the incorrect results of \cite{yeh3} and then use an example to show the faults. Then, the correct versions are established.

\vspace{2mm}

\noindent{\bf $\lq$Lemma 3' \cite{yeh3}:} Let $C$ be a {\it d-MC} candidate in $G(V,E,W)$. If there is a path between the source node and the sink node in $R(V,E,C+0(e_i))$, for $e_i\in E$, then $W(C+0(e_i))>d$.

\vspace{2mm}

\noindent{\bf $\lq$Theorem 5' \cite{yeh3}:} Let {\it C} be a {\it d-MC} candidate in {\it G(V,E,W)}. If there is a path between the source node and the sink node in $R(V, E, C + 0(e_i))$, for all $e_i \in U(C)$, then {\it C} is a real {\it d-MC}; otherwise, {\it C} is not a real {\it d-MC}.

\vspace{2mm}

A fault in $\lq$Lemma 3' is exemplified through the following example.
\begin{example}\label{ex1}
Consider Fig. 1 as a network flow. It is obvious that $C_3=\{e_1,e_3,e_4,e_6\}$ is an {\it MC} and $C_{31}=(0,2,3,1,3,3)$ is a {\it 7-MC} candidate obtained from $C_3$. It is straightforwardly seen that there is a path from node $1$ to node $4$ in $R(V,E,C_{31}+0(e_1))$ and consequently, $\lq$Lemma 3' concludes that $W(C_{31}+0(e_1))>7$, whereas $W(C_{31}+0(e_1))=6 \ngtr 7$. Hence, $\lq$Lemma 3' is incorrect.
\end{example}

Since $\lq$Theorem 5' is directly deduced from $\lq$Lemma 3', Example \ref{ex1} can be employed to illustrate the flaw in $\lq$Theorem 5' as well. In fact, $\lq$Lemma 3' and $\lq$Theorem 5' need an additional hypothesis, as $W(C)=d$. We state and prove the correct form of $\lq$Theorem 5' of \cite{yeh3} as Theorem \ref{th1} below. The correct version of $\lq$Lemma 3' can be obtained similarly by adding the assumption $W(C)=d$.

\unitlength 1.00mm 
\linethickness{1pt}
\begin{picture}(10,55)(0,35)
\put(38.75,62){\circle{5}} \put(66,73.25){\circle{5}}
\put(91.5,62){\circle{5}} \put(66,92.75){\circle{5}}
\put(38.5,62){\makebox(0,0)[cc]{1}}
\put(91.5,62){\makebox(0,0)[cc]{4}}
\put(66,92.5){\makebox(0,0)[cc]{2}}
\put(66,73){\makebox(0,0)[cc]{3}}
\put(48,80){\makebox(0,0)[cc]{$e_1$}}
\put(55,72){\makebox(0,0)[cc]{$e_2$}}
\put(65,64){\makebox(0,0)[cc]{$e_3$}}
\put(80,82){\makebox(0,0)[cc]{$e_5$}}
\put(64,84){\makebox(0,0)[cc]{$e_4$}}
\put(80,70){\makebox(0,0)[cc]{$e_6$}}
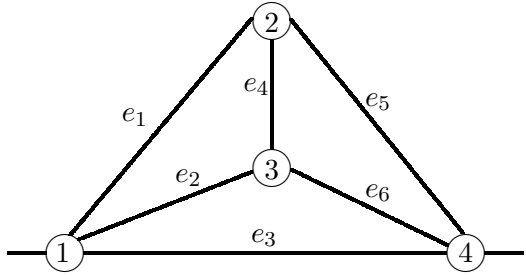
\put(65,50){\makebox(0,0)[cc]{Figure 1. The network for Example 1 with $W $=(4,2,3,1,3,3).}}

\put(41.25,62){\line(1,0){47.7}}

\put(94,62){\line(1,0){5}}

\put(36.25,62){\line(-1,0){5}}

\put(66,75.7){\line(0,1){14.5}}

\multiput(40.5,63.75)(.085576923,.033653846){267}{\line(1,0){.2}}

\multiput(68.5,73)(.069204152,-.033737024){297}{\line(1,0){.2}}

\multiput(68.5,93)(.033690658,-.042113323){673}{\line(0,-1){.2}}

\multiput(63.4,93.2)(-.0337159254,-.0405308465){710}{\line(0,-1){.2}}

\end{picture}

\begin{theorem}\label{th1}
Let {\it X} be a {\it d-MC} candidate in $G(V,E,W)$ and $W(X)=d$. {\it X} is a {\it d-MC} if and only if there is a path between the source node $1$ and the sink node $n$ in $R(V, E, (X + 0(e_i))^d)$, for all $e_i \in U(X)$.

\vspace{2mm}

\begin{pr} Suppose {\it X} is a {\it d-MC} candidate with $W(X)=d$ and there is a path between the source node $1$ and the sink node $n$ in $R(V, E, (X + 0(e_i))^d)$, for all $e_i \in U(X)$. Now, let $e$ be an arbitrary arc in $U(X)$. Since $W(X)=d$, we can send at least $d$ units of flow from node $1$ to node $n$ in $G(V,E,X + 0(e))$. Existing a path from node $1$ to node $n$ in $R(V,E,(X+ 0(e))^d)$ establishes that at least one more unit of flow can be sent from node $1$ to node $n$ as well as the $d$ units of flow previously sent in $G(V,E,X + 0(e))$. Hence, at least $d+1$ units of flow can be sent from node $1$ to node $n$, and consequently $W(X+0(e))\geq d+1>d$. Thus, according to the definition of {\it d-MC}, $X$ is a {\it d-MC}. The proof of the converse can be easily obtained.
\end{pr}
\end{theorem}

\section{On the time complexity}
Here, the complexity results of the proposed algorithm in \cite{yeh3} is meticulously investigated. Computing the correct time complexity of the proposed algorithm in \cite{yeh3}, we demonstrate that the one calculated by Yeh, $O(mnp\sigma)$ (see Theorem 6 in \cite{yeh3}), is incorrect . For convenience, the proposed algorithm in \cite{yeh3} is rewritten as Algorithm 1 below.

\begin{algorithm} An algorithm to find all the {\it d-MC}s in a limited-flow network.

Input: All {\it MC}s $C_1, C_2,..., C_p$ of a limited-flow network
$G(V,E,W)$ with the source 

node $1$ and the sink node $n$.

Output: All {\it d-MC}s.

STEP 1: Let $i=j=1$.

STEP 2: Use the Implicit Algorithm to find a feasible solution, say $X$, of the
following 

equations. If no such solution exists, then go to STEP 8.

\vspace{2mm}

$
\begin{array} {ll}
\Sigma X(e)=d,&for\ all\ e\in C_i\cr X(e)\leq W(e),&for\ all\ e\in C_i\cr X(e)=W(e),&for\ all\ e\notin C_i. 
\end{array}
$

\vspace{2mm}

STEP 3: If $X$ is the first feasible solution generated from
$C_i$ in STEP 2, then imple-

ment a {\it max-flow algorithm} to
find $W(X)$. Otherwise, employ Theorem 4 in \cite{yeh3} to

find $W(X)$.

STEP 4: If $W(X)\neq d$, then $X$ is not a {\it d-MC} candidate
and return to STEP 2 to 

find the next feasible solution.
Otherwise, $X\ $ is a$\ $ {\it d-MC} candidate, let$\ $ $C_{ij}^d\ =X,$

  $j\leftarrow j+1$ and $k=1$.

STEP 5: If $C_{ij}^d(e_k)=W(e_k)$, then go to STEP 7.

STEP 6: If there is no path from
node $1$ to node $n$ through $e_i$ in $R(V,E,C_{ij}^d+0(e_i))$, 

then $C_{ij}^d$ is not a real {\it d-MC} and return to STEP 2 to find the
next feasible solution. 

Otherwise, go to STEP 7.

STEP 7: If $k<m$, then $k\leftarrow k+1$ and go to STEP 5. Otherwise, $C_{ij}^d$ is a real {\it d-MC}.

STEP 8: If $i<p$, then $i\leftarrow i+1,\ j=1$ and go to STEP 2. Otherwise, halt.
\end{algorithm}

\subsubsection*{The correct time complexity}
In \cite{yeh3}, Yeh in Theorem 6 claimed that the time complexity of his proposed algorithm, Algorithm 1 here, was $O(mnp\sigma)$. We will show that the time complexity of Algorithm 1 is indeed $O((m^2+n^2\sqrt{m})p\sigma)$. 

It is obviously seen that the time complexity of steps 1, 4, 5, 7, and 8 is $O(1)$. In Step 2, Algorithm 1 finds a {\it d-MC} candidate using an implicit enumeration to solve the existing system of equations. Although solving the system of equations has its own complexity, the authors in [6-10] commonly disregard its complexity when calculating the time complexity of their algorithms. Therefore, we do the same in our analysis. However, the notable point is that the stated relation in Theorem 4 of \cite{yeh3} does not necessarily hold between two arbitrary {\it  d-MC} candidates obtained in Step 2. Thus, the number of usages of the {\it max-flow algorithm} in Step 3 may be more than one for each {\it MC}. In fact, we should consider the number of {\it d-MC} candidates obtained from each {\it MC} in Step 2, $p\sigma$, as the upper bound of the usage of the {\it max-flow algorithm} in Step 3 in the worst case. It should be remembered that {\it p} is the number of {\it MC}s in {\it G(V,E,W)} and the number of non-negative integer solutions generated by every {\it MC} in Step 2 of Algorithm 1 is bounded by $\sigma$. The best time complexity of the {\it max-flow algorithm} \cite{ahuja}, and so the time complexity of Step 3 is $O(n^2\sqrt{m})$. Since $ U(C_{ij}^d)$ is bounded by $m$ (the number of arcs in $E$) and the time complexity of searching for a path from the source node to the sink node is $O(m)$, the time complexity of step 6 is $O(m^2)$. Hence, the time complexity of Algorithm 1 is $O((m^2+n^2\sqrt{m})p\sigma)$ and the following theorem gives the correct result.

\begin{theorem}
The time complexity of Algorithm 1 is $O((m^2+n^2\sqrt{m})p\sigma)$.

\end{theorem}

\section{Conclusions}

Here, an existing study [Yeh WC. A new approach to the d-MC problem. Reliab Eng and Syst Saf 2002; 77(2): 201-–206.] was investigated and certain flaws in the results were pointed out. In addition, the correct versions of the results were established. Moreover, the time complexity given for the algorithm was shown to be incorrect and the correct time complexity was computed. 

\subsection*{Acknowledgements}
The authors thank the Research Council of Sharif University of Technology for its support.

\bibliographystyle{amsplain}
\bibliography{xbib}

\end{document}